\documentclass{article_Automatika}
\usepackage[utf8]{inputenc}

\usepackage{multicol}
\usepackage{multirow}
\usepackage{subfigure}
\usepackage{psfrag}
\usepackage{booktabs}
\usepackage{nicefrac}
\usepackage{amsmath}
\usepackage{amssymb}
\usepackage{gensymb}
\usepackage{graphicx}
\graphicspath{{Figures/}}

\usepackage{algorithmic}
\usepackage{algorithm}

\algsetup{
linenosize=\small,
linenodelimiter=.
}

\headertext{Projective DC Motor Control Under Disturbance Torques}{I. Zuglem, R.O. Doruk}

\author{Ismael Zuglem, R. Ozgur DORUK}
\title{Projective Control of DC Motors Under Disturbance Torques}
\titleCRO{Croatian translation of the title}

\abstract{%
In this study, we will present the design of a linear DC motor controller by projective linear qudratic servo feedback (P-LQSF) and analyze its stability through the notion of input to state stability theory. The projective control approach allows one to design an output feedback controller which approximates the eigenspectrum of a full state feedback closed loop. The performance and stability of the controllers will be analyzed both theoretically and through simulation. Apart from basic linear stability, the theoretical analysis will involve the stability of the closed loop against the disturbance torques by reflecting the closed loop as a system with the disturbance torque appearing as an input. Knowing this fact, the input-to-state stability concept is utilized as a disturbance to state stability approach and the designs are analyzed accordingly. The overall products are demonstrated by MATLAB based simulations.      
}%


\begin{document}

\maketitle

\section{Introduction}\label{sec:Intro}
Direct Current Machines are widely used torque transducers in mechanical systems the applications of which range from automotive, robotics, pneumatic and hydraulic systems and various biomedical engineering applications. The simplest version of a DC motor involves a permanent magnet rotor and an armature winding (stator) which is often the case when one has a brushed DC motor \cite{beaty1998electric,krishnan2001electric}. These can be modeled according to the fundamental circuit theories and often available in control systems textbooks such as \cite{Ogata:2001:MCE:516039}. The brushless DC motors which have a permanent magnet stator but a wound rotor is actually an AC motor. It will require a dedicated driver circuitry to be operated from a DC supply \cite{pillay1989modeling}. The dynamical characteristics of brushed or brushless DC motors are similar \cite{rao2012mathematical}. A proper and beneficial DC motor application will require a position or speed controller so that the desired performances are obtained (constant or tracked speed and/or position). They are also a well established class of mechanical systems suited for control system development and numerous researches are available in literature. 

Concerning the control approaches one can note that regardless of targeting position and speed most of the motor controller designs involves proportional+integrator+derivative group (PID or PI) of controllers. Some related examples can be found in \cite{alexandridis2014modified,huang2008pc,kanojiya2012optimal,lin1994comparison,liu2009real,meshram2012tuning,rao1996neuro,yu2004optimal,zhou2008dc}. The PID group of controllers are structured control laws that can be tuned according to various methodologies such as Zeiger-Nichols charts \cite{meshram2012tuning}, optimization \cite{kanojiya2012optimal,yu2004optimal} and even neuroadaptive \cite{rao1996neuro} techniques. Some other control related studies in the literature are about fuzzy logic based motor controls \cite{lin1994comparison}, a Kalman Filter based example \cite{praesomboon2009sensorless} and another application using optimal state feedback \cite{ruderman2008optimal}.

Almost all control approaches need to implement a feedback from all or part of its state variables (position, angular velocity/speed, armature current etc.). Depending on the application some of the state variables may or may not be available for measurement. Factors affecting this availability may be the cost, the feasibility of the usage of certain instruments such as tachometers, encoders, current or torque sensors. In the literature, there are sensorless control approaches such as \cite{afjei2007sensorless,montanari2000speed,praesomboon2009sensorless}. These aim at the control of motor dynamics without the employment of a position or speed sensor. Such approaches generally require the utilization of an observer \cite{Ogata:2001:MCE:516039} such as a Kalman Filter \cite{kalman1960new,julier1997new,van2001square,wan2000unscented}. Elimination of an observer/filter means employment of a static output feedback approach which is lack of a profound systematic knowledge. However thanks to \cite{medanic1983design,medanic1985design,medanic1988design}, the flexibility of a full state feedback control can be reflected (or projected) to an output feedback by a simple orthogonal projection operation from the state space of the full state feedback to the state space of the output only feedback. This approach is formerly used in aerospace applications \cite{wise1991approximating,wise1992optimal}, some process control applications \cite{saif1989suboptimal,srinivasan1994projective} and also as a Dynamic PI Control tuning helper \cite{meo1986design}.

In this research, we will utilize the approach presented in \cite{medanic1983design} to dc motor control assuming that the feedback from the armature current is not available. This is a practically possible situation as the measurement of a current through a sense resistor followed by an signal amplification might lead to accumulation of unwanted noise. We will present a position and a speed control application which have a feedback only from the speed/position tracking errors. 

Motors are often subject to disturbance torques when operated at harsh environments such as non-uniform fluid flows over the propeller or bad lubrication of the bearings etc. The level of disturbance torques might be a threat to the stability of the closed loop control system. The last part of this research is to deal with the effects of the available disturbance torques. The disturbance torques can be considered as an exogenous input to the closed loop system and thus the notion of input-to-state stability \cite{sontag1995characterizations,sontag1995input} can be considered as a useful analysis approach. The disturbance decoupling concept \cite{willems1981disturbance,weiland1989almost,saberi1987output} is an extension of input-to-state stability property which analyses the disturbance quenching capability of the closed loop control system. In this research, we will try to assess the conditions of input-to-state stability treating the disturbances as inputs to the closed loop. 

The main contribution of this research to the literature can be summarized as follows:
\begin{itemize}
\item Application of \emph{Projective Control Approach} to the design of electric motor control systems
\item Analysis of the disturbance handling capabilities of a closed loop motor control system through the utilization of input-to-state stability. 
\end{itemize} 
The demonstration of the results will be performed by numerical simulation of the designs. MATLAB is the main computational environment in this study.               
\section{Direct Current Motor Model}\label{sec:DCMotor}
The direct current electrical machines generally involves the dynamics of the shaft position $ \theta(t) $ (in radians or degrees) and angular velocity $ \omega(t) $ (in $ \nicefrac{rad}{sec} $) and armature current $ i_a(t) $ (in Amperes). We are not talking about the field winding here.
\begin{subequations}\label{eq:DCM_EQ}
\begin{align} 
\dot{\theta}&=\omega\\ 
\dot{\omega}&=\alpha\omega+\beta i_a+\nu \tau_L \\
\dot{i_a}&=\gamma\omega+\rho i_a+sV_a
\end{align} 
\end{subequations}

\noindent where $ \alpha=-\frac{B}{J} $, $ \beta=\frac{K_i}{J} $, $ \gamma=-\frac{K_b}{L_a} $, $ \rho=-\frac{R_a}{L_a} $, $ s=\frac{1}{La} $ and $ \nu=\frac{1}{J} $. In \textbf{Table \ref{tab:DCMPARAM}}, one can see the definitions and  nominal values of a particular DC motor. In state space form one will have the following equations:
\begin{subequations} \label{SSDC}
\begin{equation} \label{SSDC_POSITION}
\begin{bmatrix}
\dot{\theta}\\ \dot{\omega} \\ \dot{i_a}
\end{bmatrix}=\begin{bmatrix}
0 & 1 & 0 \\
0 & \alpha & \beta \\
0 & \gamma & \rho 
\end{bmatrix}\begin{bmatrix}
\theta\\ \omega \\ i_a
\end{bmatrix}+\begin{bmatrix}
0 \\ 0 \\ s
\end{bmatrix}V_a+\begin{bmatrix}
0 \\ \nu \\ 0
\end{bmatrix}\tau_L
\end{equation}
\begin{equation} \label{SSDC_SPEED}
\begin{bmatrix}
\dot{\omega} \\ \dot{i_a}
\end{bmatrix}=\begin{bmatrix}
\alpha & \beta \\
\gamma & \rho 
\end{bmatrix}\begin{bmatrix}
\omega \\ i_a
\end{bmatrix}+\begin{bmatrix}
0 \\ s
\end{bmatrix}V_a+\begin{bmatrix}
\nu \\ 0
\end{bmatrix}\tau_L
\end{equation}
\end{subequations}
The equation in \eqref{SSDC} has two subsections. The state space representation in \eqref{SSDC_POSITION} is intended for position control where as \eqref{SSDC_SPEED} is intended for speed control.   
\begin{table}[!htb]
\caption{DC Motor Parameters and Definitions}
\centering
\begin{tabular}{|c|c|c|}
\hline
Definition & Symbol & Value \\ \hline
Inertia of the Load & $J$ & $ 0.01\ kg\cdot\mathrm{m}^2 $  \\ \hline
Viscous friction coefficient & $B$ & $ 0.1\ \mathrm{N\cdot m\cdot \nicefrac{sec}{rad}} $ \\ \hline
Armature Resistance & $R_a$ & $ 1\ \Omega $   \\ \hline
Armature Inductance & $L_a$ & $ 0.5\ \mathrm{H} $   \\ \hline
Torque Constant & $K_i$ & $ 0.01\ \mathrm{N\cdot \nicefrac{m}{A}} $   \\ \hline
Back-EMF Constant & $K_b$ & $ 0.01\ \mathrm{V\cdot \nicefrac{sec}{rad}} $   \\ \hline
\end{tabular}
\label{tab:DCMPARAM}
\end{table} 

\noindent The term $ \tau_L $ stands for the load or disturbance torques due to certain restrictive factors. These may be due to aerodynamic factors in a propeller or an unpredictable friction on the shaft etc. The closed loop of the motor model in \eqref{eq:DCM_EQ} will still have the disturbance torque $ \tau_L $ as input. Because of that we will make use of the input-to-state stability approaches in [reference] to analyze our motor controller under the disturbance torques. 

Considering the control approaches there are various alternatives as stated in \textbf{Section \ref{sec:Intro}}. The main issue about the full state feedback techniques is that some state variables can not be measured. In these cases either an observer [luenberger and kalman filters] should be used. In DC motor models such as \eqref{SSDC}, the armature current which determines the torque through the relation $ \tau=K_i\times i_a $ may not be easy to sense continuously. Though devices such as low resistance sense resistors are often used in current measurement, their utilization in control requires amplification (such as OP-AMPS or Instrumentation-Amplifiers) which may bring noise and offset adjustment requirements. Thus, a control approach that does not need current feedback may benefit from being free of those issues. Apart from these, lower number of instruments will be required which is a cost reduction measure. 

In the next section we will introduce a methodology called as projective control which can be applied to design an output feedback controller. 

Classical control techniques based on transfer function and compensation approaches [reference] will not be considered here as they are well established classical methodologies.

\section{LINEAR QUADRATIC PROJECTIVE CONTROL APPROACH} \label{LQPC}
Consider the following linear model in state space form:
\begin{equation}
\dot{x}=Ax+Bu \label{LINSSGEN}
\end{equation}
where $ x\in\Re^n $, $ u\in\Re $, $ A\in\Re^{n\times n} $ and $ B\in\Re^n $. The linear quadratic regulation is to find a control law $ u $ which will minimize the following quadratic performance index:
\begin{equation}
J=\int\limits_{0}^{\infty}(x^TQx+u^TRu)dt \label{QUADPI}
\end{equation}
along the trajectories of \eqref{LINSSGEN}. The control law is generally in a \emph{full state feedback form} of $ u=-Kx $ where $ K=R^{-1}B^TP $. The matrix $ P $ is a positive definite symmetric solution of the algebraic Riccati equation which is:
\begin{equation}
A^T P + P A - (P B) R^{-1} (B^T P) + Q = 0 \label{LQRICCAT}
\end{equation}    
where $ Q $ is also a symmetric positive definite matrix. In many applications, this matrix can be chosen as $ Q=qI_{n\times n} $ where $ q $ is a positive constant. The closed loop dynamics when the feedback $ u=-Kx $ is applied to \eqref{LINSSGEN} will be obtained as:
\begin{equation}
\dot{x}=(A-BK)x \label{SSGENCLOSED}
\end{equation}
Provided that, the pair $ (A,B) $ is complete state controllable and a positive definite symmetric solution $ P $ is found for the Riccati equation in \eqref{LQRICCAT}, the eigenvalues of \eqref{SSGENCLOSED} should have negative real parts. 

In many applications, feedback from all state variables are not possible. At least one or two variables are not physically measurable. This is often a burden in employment of state space based modern control techniques because one needs to employ an output feedback. In some cases these issues are handled by utilization of transfer function based techniques but these approaches hide the benefits obtained from state space based techniques. Because of that researchers worked on approaches which allows output feedback directly from state space representations. One of the related works is [reference] which makes use of the orthogonal projection theorem in mathematics to project the closed loop eigenspectrum of a full state feedback control approach to the closed loop eigenspectrum of an output feedback just from the measurable states. To further develop this approach one can rewrite \eqref{LINSSGEN} with the output to be fed back as follows:
\begin{subequations} \label{LINSS_OUTPUT}
\begin{equation} 
\label{LINSS_STATE}
\dot{x}=Ax+Bu
\end{equation}
\begin{equation} 
\label{LINSS_MEASOUT}
y=Cx
\end{equation}
\end{subequations}
where $ C $ in \eqref{LINSS_MEASOUT} is a matrix which filters out the measurable states to the output. It might be formed from ones and zeros in order to create direct feedback from the measurable states in $ x $. Sometimes even all states in $ x $ are measurable (full state feedback is practically applicable) the design might require feedback from a combination of the state variables. In such cases projective control approach is also an applicable method.    

The output feedback from $ y $ will be applied as $ u=-K_oCy $ where $ K_o\in\Re^{n\times r} $ is the output feedback gain and $ C\in\Re^{r\times n} $. In those relationships, $ r $ is the number of measurable states (or outputs in general if some variables in $ y $ are linear combination of the states). The closed loop will have the following dynamics:
\begin{equation}
\dot{x}=(A-BK_oC)x \label{SSOUTFBCLOSED}
\end{equation}  
The above may or may not have eigenvalues with negative real parts. However, it is well known from [reference] that $ r $ number of the eigenvalues can be guaranteed to be made stable whereas the rest $ n-r $ eigenvalues are not manipulatable. Of course, this fact does not mean that the output feedback always produces unstable designs. 

The projective control approach is developed from the eigenspectrum relations between \eqref{SSOUTFBCLOSED} and \eqref{SSGENCLOSED}. One can write the eigendecomposition relation for the full state feedback \eqref{SSGENCLOSED} as:
\begin{equation} \label{SSFFBEIGENSPEC}
(A-BK)V=V\Lambda
\end{equation}
and for the output feedback closed loop \eqref{SSOUTFBCLOSED} as:
\begin{equation} \label{SSOUFBEIGENSPEC}
(A-BK_oC)V_r=V_r\Lambda_r 
\end{equation}
where $ V\in\Re^{n\times n} $ is the eigenvector matrix consisting of the eigenvalues of \eqref{SSGENCLOSED} denoted by $ \Lambda\in\Re^{n\times n} $. In \eqref{SSOUFBEIGENSPEC}, $ \Lambda_r\in\Re^{r\times r} $ denotes the $ r $ eigenvalues chosen from $ \Lambda $ that are to be retained when the feedback $ u=-K_oy $ is applied. The matrix $ V_r\in\Re^{n\times r} $ consists of eigenvectors corresponding to the eigenvalues in $ \Lambda_r $. Note that, as $ \Lambda_r $ is an $ r $-element subset of $ \Lambda $ one can also write the following relation:
\begin{equation} \label{SSFFBEIGENSPEC_PROJEIG}
(A-BK)V_r=V_r\Lambda_r
\end{equation}
and the above can be equated to \eqref{SSOUFBEIGENSPEC} and thus:
\begin{equation} \label{PROJCONTROL_DERIVTN}
(A-BK)V_r=V_r\Lambda_r=(A-BK_oC)V_r
\end{equation}  
As a result, one will be able to write the relationship between $ K_o $ and $ K $ as shown below:
\begin{equation} \label{PROJCONTROL_EQUATION}
K_o=KV_r(CV_r)^{-1}
\end{equation}    
So one can first find a full state feedback control law as $ u=-Kx $ from \eqref{LQRICCAT} and then by applying \eqref{PROJCONTROL_EQUATION} to the full state feedback gain $ K $. 

Selection of the eigenvalues to be retained $ (\Lambda_r) $ among the full state feedback closed loop eigenvalues $ (\Lambda) $ depends on the number of available feedback lines $ (r) $, the nature of the poles i.e. whether they are complex conjugates or real and they are dominance. First of all, the r dominant eigenvalues among $ \Lambda $ should be considered. If the number complex eigenvalues restricts this choice then other eigenvalues may also be selected. 
  
\section{CONTROL OF DIRECT CURRENT MOTORS} 
Position and speed control of direct current motors may require different configurations. For example position dynamics in \eqref{SSDC_POSITION} involve a natural integration which helps in the elimination of any steady state error. The speed (angular velocity) dynamics \eqref{SSDC_SPEED} on the contrary does not have any natural integrator which may lead to a steady state error. In order to overcome this issue one will need to add an artificial integrator to \eqref{SSDC_SPEED}. In this section, we will present the speed and position control of the DC motors with the aid of projective control approach of \textbf{Section \ref{LQPC}}. 
\subsection{Speed Control of DC Motors by Projective Control}  \label{SPD_CONTROL_SEC}
As the speed dynamics \eqref{SSDC_SPEED} have no natural integrator one has to add an artificial integrator. This can be performed by adding this integrator to the forward path. This means that one has to integrate the tracking error as shown below:
\begin{equation}
\dot{\epsilon}=\omega-\omega_r
\end{equation} 
where $ \omega_r $ is the reference (desired) angular velocity and $ \omega $ is the actual speed (angular velocity) of the motor shaft. So the state space representation of the motor model for speed control will be obtained from the combination of the above and \eqref{SSDC_SPEED} as:
\begin{equation} \label{SPDCONTROL_SS}
\begin{bmatrix}
\dot{\epsilon}\\ \dot{\omega} \\ \dot{i_a}
\end{bmatrix}=\begin{bmatrix}
0 & 1 & 0 \\
0 & \alpha & \beta \\
0 & \gamma & \rho 
\end{bmatrix}\begin{bmatrix}
\epsilon\\ \omega \\ i_a
\end{bmatrix}+\begin{bmatrix}
0 \\ 0 \\ s
\end{bmatrix}V_a+\begin{bmatrix}
0 \\ \nu \\ 0
\end{bmatrix}\tau_L+\begin{bmatrix}
-1 \\ 0 \\ 0
\end{bmatrix}\omega_r
\end{equation}  
as stated in [reference ogata olabilir] if the reference speed is slowly changing or constant the difference between the state variables $ \epsilon(t) $, $ \omega(t) $, $ i_a(t) $ and their steady state values $ \epsilon(\infty) $, $ \omega(\infty) $, $ i_a(\infty) $ will not involve reference speed $ \omega_r $. So one can write the following:
\begin{equation}
\begin{bmatrix} \label{SPDCONTROL_SS_NO_TAU}
\dot{e}_\epsilon\\ \dot{e}_\omega \\ \dot{e}_{i_a}
\end{bmatrix}=\begin{bmatrix}
0 & 1 & 0 \\
0 & \alpha & \beta \\
0 & \gamma & \rho 
\end{bmatrix}\begin{bmatrix}
e_\epsilon\\ e_\omega \\ e_{i_a}
\end{bmatrix}+\begin{bmatrix}
0 \\ 0 \\ s
\end{bmatrix}V_a
\end{equation}
where $ e_\epsilon=\epsilon(t)-\epsilon(\infty) $, $ e_\omega=\omega(t)-\omega(\infty) $, $ e_{i_a}=i_a(t)-i_a(\infty) $. The control law can be rewritten as:
\begin{equation}
V_a=-K^fx=-\begin{bmatrix} \label{FSFB_DC_MOTOR_SPD}
k^f_\epsilon & k^f_\omega & k^f_{i_a}
\end{bmatrix}\begin{bmatrix}
e_\epsilon \\ e_\omega \\ e_{i_a}
\end{bmatrix}
\end{equation} 
for the case of full state feedback and,
\begin{equation}
V_a=-K^oCx=-\begin{bmatrix} \label{OUFB_DC_MOTOR_SPD}
k^o_\epsilon & k^o_\omega
\end{bmatrix}\begin{bmatrix}
e_\epsilon \\ e_\omega
\end{bmatrix}
\end{equation} 
for output feedback control. In \eqref{FSFB_DC_MOTOR_SPD} and \eqref{OUFB_DC_MOTOR_SPD} $, e=\begin{bmatrix}e_\epsilon & e_\omega & e_{i_a}\end{bmatrix}^T $. In \eqref{OUFB_DC_MOTOR_SPD} the output matrix $ C $ can be written as:
\begin{equation}
C=\begin{bmatrix} \label{CMATR}
1 & 0 & 0\\
0 & 1 & 0
\end{bmatrix}
\end{equation}
The above means that, we are having a feedback from $ e_\epsilon $ and $ e_\omega $ which are the state variables related to motor shaft speed $ (\omega) $.  
Note that, in the above discussion the load/disturbance torque $ \tau_L $ is currently not taken into account and left for the disturbance analysis section.

The closed loop dynamics when full state feedback in \eqref{FSFB_DC_MOTOR_SPD} is employed, will appear as:
\begin{subequations} \label{CL_FS_DYN}
\begin{equation} \label{CL_FS_DYN_LTR}
\dot{e}=(A-BK^f)e
\end{equation}
\begin{equation} \label{CL_FS_MATRIX}
A-BK^f=\begin{bmatrix}
0 & 1 & 0\\
0 & \alpha & \beta\\
-sk^f_\epsilon & \gamma-sk^f_\omega & \rho-sk^f_{i_a}
\end{bmatrix}
\end{equation}
\end{subequations}
The next step is to apply projective control equation in \eqref{PROJCONTROL_EQUATION}. In the speed control of DC motors by projective control approach, as understood from \eqref{CMATR} one will have only two feedback variables. This means that we can only retain two eigenvalues from the closed loop spectrum of \eqref{CL_FS_MATRIX}. Due to the odd size of \eqref{CL_FS_MATRIX}, if it has two complex eigenvalues one has to select them as the retained eigenvalues. If all the eigenvalues are real the two dominant ones should be preferred as retained eigenvalues. After the application of \eqref{PROJCONTROL_EQUATION} one should check the closed loop eigenvalues of the output feedback using the following: 
\begin{subequations} \label{CL_OU_DYN}
\begin{equation} \label{CL_OU_DYN_LTR}
\dot{e}=(A-BK^oC)e
\end{equation}  
\begin{equation} \label{CL_OU_MATRIX}
A-BK^oC=\begin{bmatrix}
0 & 1 & 0\\
0 & \alpha & \beta\\
-sk^f_\epsilon & \gamma-sk^f_\omega & \rho
\end{bmatrix}
\end{equation}
\end{subequations}

The information up to this point will be applied in \textbf{Section eklenecek} to a numerical application.  
\subsection{Position Control of DC Motors by Projective Controls} 
The position dynamics in \eqref{SSDC_POSITION} has a natural integrator as position $ \theta(t) $ is the integration of the angular velocity $ \omega(t) $. So one does not have to add any sort of artificial integrators. 
\begin{equation} \label{PSNCONTROL_SS}
\begin{bmatrix}
\dot{\theta}\\ \dot{\omega} \\ \dot{i_a}
\end{bmatrix}=\begin{bmatrix}
0 & 1 & 0 \\
0 & \alpha & \beta \\
0 & \gamma & \rho 
\end{bmatrix}\begin{bmatrix}
\theta\\ \omega \\ i_a
\end{bmatrix}+\begin{bmatrix}
0 \\ 0 \\ s
\end{bmatrix}V_a+\begin{bmatrix}
0 \\ \nu \\ 0
\end{bmatrix}\tau_L+\begin{bmatrix}
-1 \\ 0 \\ 0
\end{bmatrix}\theta_r
\end{equation}  
where $ \theta_r $ is the reference or desired position. Other variables are exactly the same as that for \eqref{SPDCONTROL_SS}. It should also be noted that, \eqref{PSNCONTROL_SS} and \eqref{SPDCONTROL_SS} are exactly the same except the error integral variable $ \epsilon $ which is replaced by position $ \theta $ and the reference speed variable $ \omega_r $ that is replaced by $ \theta_r $.
Because of that, the details given in \textbf{Section \ref{SPD_CONTROL_SEC}} can be directly applied here provided that the designer is careful about the state variables. The error dynamics in \eqref{SPDCONTROL_SS_NO_TAU} is replaced by:

\begin{equation}
\begin{bmatrix} \label{PSNCONTROL_SS_NO_TAU}
\dot{e}_\theta\\ \dot{e}_\omega \\ \dot{e}_{i_a}
\end{bmatrix}=\begin{bmatrix}
0 & 1 & 0 \\
0 & \alpha & \beta \\
0 & \gamma & \rho 
\end{bmatrix}\begin{bmatrix}
e_\theta\\ e_\omega \\ e_{i_a}
\end{bmatrix}+\begin{bmatrix}
0 \\ 0 \\ s
\end{bmatrix}V_a
\end{equation}
and the control equations \eqref{FSFB_DC_MOTOR_SPD} and \eqref{OUFB_DC_MOTOR_SPD} will be:
\begin{equation}
V_a=-K^fe=-\begin{bmatrix} \label{FSFB_DC_MOTOR_PSN}
k^f_\theta & k^f_\omega & k^f_{i_a}
\end{bmatrix}\begin{bmatrix}
e_\theta \\ e_\omega \\ e_{i_a}
\end{bmatrix}
\end{equation} 
for the case of full state feedback. As one should have the feedback from position tracking error $ e_\theta $ and the speed error w.r.to the steady state $ e_\omega $ the feedback matrix should be same as that of \eqref{CMATR} (first two elements of the state vector).    
\begin{equation}
V_a=-K^oCe=-\begin{bmatrix} \label{OUFB_DC_MOTOR_PSN}
k^o_\theta & k^o_\omega
\end{bmatrix}\begin{bmatrix}
e_\theta \\ e_\omega
\end{bmatrix}
\end{equation}
Finally, the closed loop dynamics \eqref{CL_FS_DYN} and \eqref{CL_OU_DYN} will become:
\begin{subequations} \label{CL_FS_DYN_PSN}
\begin{equation} \label{CL_FS_DYN_LTR_PSN}
\dot{e}=(A-BK^f)e
\end{equation}
\begin{equation} \label{CL_FS_MATRIX_PSN}
A-BK^f=\begin{bmatrix}
0 & 1 & 0\\
0 & \alpha & \beta\\
-sk^f_\theta & \gamma-sk^f_\omega & \rho-sk^f_{i_a}
\end{bmatrix}
\end{equation}
\end{subequations}
And for the output feedback closed loop will be:
\begin{subequations} \label{CL_OU_DYN_PSN}
\begin{equation} \label{CL_OU_DYN_LTR_PSN}
\dot{e}=(A-BK^oC)e
\end{equation}
\begin{equation} \label{CL_OU_MATRIX_PSN}
A-BK^oC=\begin{bmatrix}
0 & 1 & 0\\
0 & \alpha & \beta\\
-sk^f_\theta & \gamma-sk^f_\omega & \rho
\end{bmatrix}
\end{equation}
\end{subequations}
One will need to verify the closed loop's stability through \eqref{CL_FS_MATRIX_PSN}. The selection of the retained eigenvalues are subject to the same rules as described in the speed control. A numerical application will be given in \textbf{Section eklenecek}.
\section{INPUT-TO-STATE STABILITY} 
\subsection{Theory} \label{THEORY_INPSTAB}
In this section, we will introduce the input-to-state stability concept and present its applicability in the analysis of the overall stability of the closed loop motor controller against disturbance torques. Before proceeding it will be beneficial to give some definitions [reference gerekli/eduardo sonntag ve bazi diger yararli paperlar olabilir definitionlarin ve theoremlerin icindede verilmeli]. 
\begin{definition}[Class $ \mathcal{K} $ Functions]
These are a class of all functions $ \eta:\Re_+\rightarrow\Re_+ $ satisfying the following conditions:
\begin{enumerate}
\item $ \eta(0)=0 $
\item $ \eta(.) $ is continuous
\item $ \eta(.) $ is strictly increasing
\end{enumerate}
\end{definition}
\begin{definition}[Class $ \mathcal{K}_\infty $ Functions]
$ \xi(p) $ will be of class $ \mathcal{K}_\infty $ iff $ \xi(.) $ is of class $ \mathcal{K} $ and $ \xi(p)\rightarrow\infty $ when $ p\rightarrow\infty $.   
\end{definition}
\begin{definition}[Storage Functions] \label{STORAGEDEF}
The function $ W(x):\Re^n\rightarrow\Re_+ $ with $ x\in\Re^n $is said to be a storage Lyapunov function if it satisfies the following conditions:
\begin{enumerate}
\item $ W $ is continuously differentiable
\item $ W $ is radially unbounded i.e. $ W(x)\rightarrow\infty $ when $ x\rightarrow\infty $
\item $ W $ is a positive definite function i.e. $ W(0)=0 $ and $ W(x)>0 $ when $ x\ne0
 $.
\end{enumerate}
\end{definition}
\begin{theorem}[Stability of an autonomous system] \label{LYAPSTABTHM}
Suppose that a system defined by the following differential equation:
\begin{equation}
\dot{x}=f(x) \label{NLINSYSTEM}
\end{equation}
and also suppose that $ f(0)=0 $. The equilibrium point $ x=0 $ will be stable in the sense of Lyapunov, if a function $ W(x) $ satisfying the properties given in \textbf{Definition \ref{STORAGEDEF}} and:
\begin{equation}
\frac{\partial W(x)}{\partial x}f(x)\leq -\eta(|x|)
\end{equation} 
where $ \eta(|x|) $ is a class $ \mathcal{K}_\infty $ function.    
\end{theorem}
\textbf{Theorem \ref{LYAPSTABTHM}} is valid when the system has no exogenous inputs. When the system in \eqref{NLINSYSTEM} has exogenous or normal inputs ($ u $) one needs to take it into consideration. 
In this case we will need to define an \emph{ISS-Lyapunov} function. See the definition below:
\begin{definition}[\emph{ISS-Lyapunov} Functions] \label{ISSLYAPDEF}
An ISS - Lyapunov function $ W(x) $ is a type of storage function which satisfies the properties given in \textbf{Definition \ref{STORAGEDEF}} and the one shown below:
\begin{equation} \label{DISSIPISSINEQ}
\frac{\partial W(x)}{\partial x}f(x,u)\leq -\eta(|x|)+\theta(|u|)
\end{equation}
where $ \eta(.) $ and $ \theta(.) $ are class $ \mathcal{K}_\infty $ functions.  
\end{definition} 
\begin{theorem}[Input-to-State Stability (ISS)]
\label{ISSTHEOREM}
A system defined by the following general differential equation:
\begin{equation}
\dot{x}=f(x,u) \label{NLINSYSTEMU}
\end{equation} 
will be input to state stable (ISS) if there exist a $ W(x) $ satisfying the properties in \textbf{Definitions \ref{STORAGEDEF} and \ref{ISSLYAPDEF}} for the system in \eqref{NLINSYSTEMU}. 
\end{theorem}
In general, the above case can be associated with the dissipation concept as \eqref{DISSIPISSINEQ} is a dissipation inequality with the storage function $ W(x) $ and the supply function $ \sigma(x,u)=\theta(|u|)-\eta(|x|) $. 

Input-to-state stability can be viewed as a disturbance-to-state stability when the closed loop formed by applying a feedback of the form $ u=-k(x) $ to a general system with a exogenous disturbance input $ n(t) $:
\begin{equation}
\dot{x}=f(x,u,n) \label{NLINSYSTEMUW}
\end{equation}
and when the feedback is applied the above system will become $ \dot{x}=f(x,-k(x),n) $. Thus, the closed loop can be thought of a system with input $ n $. Then the input to state stability condition in \eqref{DISSIPISSINEQ} can be rewritten as:
\begin{equation}
\frac{\partial W(x)}{\partial x}f(x,w)\leq -\eta(|x|)+\theta(|w|) \label{DISTURBTOSTAT}
\end{equation} 
with the definitions of $ \alpha $, $ \theta $ are same as that of \eqref{DISSIPISSINEQ}. The above condition is called as Disturbance-to-State condition.   
\subsection{Disturbance-to-State Stability for DC Motor Control}
In this section, one will be able to see how the theory developed in \textbf{Section \ref{THEORY_INPSTAB}} is applied to the analysis of stability against the disturbance torques exerted on the DC motor shaft and load. To achieve this goal one should first take the closed loop dynamics in \eqref{CL_OU_DYN} or \eqref{CL_OU_DYN_PSN}. However, referring to equation \eqref{SPDCONTROL_SS} or \eqref{PSNCONTROL_SS} one should modify the closed loop dynamics to include the disturbance torque $ \tau_L $. That is:
\begin{equation}
\label{OU_FB_DISTURBANCE_ABG}
\dot{e}=(A-BK^oC)e+G\tau_L
\end{equation} 
where $ G=\begin{bmatrix}0 & -1 & 0\end{bmatrix}^T $. Before continuing some additional information should be presented. 
\begin{definition}[Quadratic Forms] \label{QUADFORMDEF}
For any symmetric matrix $ P\in\Re^{n\times n} $, the form $ x^TPx $ will be called as a quadratic form. 
\end{definition}
\begin{theorem}[Lower and Upper Bounds] 
\label{QFORMBOUNDS} 
For any quad-ratic form defined in \textbf{Definition \ref{QUADFORMDEF}} one can define the following lower and upper bounds:
\begin{equation}
\label{EQQUADBOUND}
\lambda_{\min}(P)\leq x^TPx \leq \lambda_{\max}(P)
\end{equation} 
where $ \lambda_{\min}(P) $ and $ \lambda_{\max}(P) $ are the minimum and maximum eigenvalues of the matrix $ P $ respectively. 
\end{theorem}
Now considering the storage function concept in \textbf{Section \ref{THEORY_INPSTAB}} one can define the following as a quadratic storage function:
\begin{equation}
W(e)=\frac{1}{2}e^Te \label{LYAPFORe}
\end{equation} 
and one can also write the rate of change of $ W(e) $ along the trajectories of $ e $ as:
\begin{equation}
\label{RATEOFCHANGEOFW}
\dot{W}(e)=\frac{\partial W(e)}{\partial e}\dot{e}=e^T\dot{e}
\end{equation}
and substituting from \eqref{OU_FB_DISTURBANCE_ABG} one will be able to obtain:
\begin{multline}
\label{WDOTOPENFORM}
\dot{W}(e)=e^T\left\lbrace (A-BK^oC)e+G\tau_L \right\rbrace\\= e^T(A-BK^oC)e+e^TG\tau_L
\end{multline}
Using \textbf{Theorem \ref{QFORMBOUNDS}}, one can convert the above to an inequality as:
\begin{equation}
\label{WDOTOPENFORM_INEQ}
\dot{W}(e)\leq\lambda_{\max}(A-BK^oC)\left\|e \right\|^2+e^TG\tau_L 
\end{equation}
In order to go further, we will need to deal with the term $ e^TG\tau_L $. It is pretty obvious that $ (e-G\tau_L)^T(e-G\tau_L)\geq0 $. We can expand this term as shown below:
\begin{multline}
\label{SECOND_TERM_1}
(e-G\tau_L)^T(e-G\tau_L)=\\e^Te-e^TG\tau_L-\tau_L^TG^Te+\tau_L^TG^TG\tau_L\geq0
\end{multline} 
Compiling the right side of the equation:
\begin{equation}
\label{SECOND_TERM_2}
e^Te+\tau_L^TG^TG\tau_L \geq e^TG\tau_L+\tau_L^TG^Te
\end{equation}
when $ \tau_L $ is a scalar as in \eqref{SSDC} $ e^TG\tau_L=\tau_L^TG^Te $ so one can write the following:
\begin{equation}
\label{SECOND_TERM_3}
e^Te+\tau_L^TG^TG\tau_L \geq 2e^TG\tau_L
\end{equation}
and
\begin{equation}
\label{SECOND_TERM_4}
\frac{1}{2}(e^Te+\tau_L^TG^TG\tau_L) \geq e^TG\tau_L
\end{equation}
Using \eqref{SECOND_TERM_4} one can rewrite \eqref{WDOTOPENFORM_INEQ} as:
\begin{multline}
\dot{W}(e)\\\leq\lambda_{\max}(A-BK^oC)\left\|e \right\|^2+\frac{1}{2}e^Te+\frac{1}{2}\tau_L^TG^TG\tau_L
\end{multline}
As it is known that $ e^Te=\left\| e \right\|^2  $ in the sense of $ \mathcal{L}_2 $ norms, the above inequality can be rewritten as:
\begin{multline}
\dot{W}(e)\\\leq\lambda_{\max}(A-BK^oC)\left\|e \right\|^2+\frac{1}{2}\left\| e \right\|^2+\frac{1}{2}\tau_L^TG^TG\tau_L
\end{multline} 
The first two terms of the above can be combined as:
\begin{multline}
\dot{W}(e)\\\leq\left[ \frac{1}{2}+\lambda_{\max}(A-BK^oC)\right] \left\|e \right\|^2+\frac{1}{2}\tau_L^TG^TG\tau_L
\end{multline}
Finally using \textbf{Definition \ref{QFORMBOUNDS}} once again, the above inequality is finalized as:
\begin{multline}
\label{DIST_TO_STAT_INEQ}
\dot{W}(e)\\\leq\left[ \frac{1}{2}+\lambda_{\max}(A-BK^oC)\right] \left\|e \right\|^2+\frac{1}{2}\lambda_{\max}(G^TG)\left\|\tau_L \right\|^2
\end{multline} 
Looking at \eqref{DIST_TO_STAT_INEQ}, one will easily note that it resembles \eqref{DISTURBTOSTAT}. The input-to-state stability theorem (\textbf{Theorem \ref{ISSTHEOREM}}) dictates that $ \frac{1}{2}+\lambda_{\max}(A-BK^oC)<0 $ and $ \lambda_{\max}(G^TG)>0 $. One can now state the following theorem:
\begin{theorem}[Stability Against Disturbances]
\label{DSSTHEOREM}
The closed loop controlled DC motor modeled by equations \eqref{SPDCONTROL_SS} or \eqref{PSNCONTROL_SS} with the control law presented in \eqref{OUFB_DC_MOTOR_SPD} or \eqref{OUFB_DC_MOTOR_PSN} will be disturbance-to-state $( \tau_L $-to-$ e )$ stable if the following conditions are satisfied:
\begin{subequations}
\begin{equation}
\label{DSS_CONDITA}
\lambda_{\max}(A-BK^oC)<-\frac{1}{2}
\end{equation}
\begin{equation}
\label{DSS_CONDITB}
\lambda_{\max}(G^TG)>0
\end{equation}
\end{subequations}
\end{theorem}  
One should note that, \textbf{Theorem \ref{DSSTHEOREM}} is a sufficient condition not a necessary one. The $ G^TG $ matrix is evaluated as 
\begin{equation}
G^TG=\nu^2
\end{equation}
which is a scalar and $ \lambda_{\max}(G^TG)=\nu^2>0 $. So condition \eqref{DSS_CONDITB} is always satisfied. The condition in \eqref{DSS_CONDITA} requires numerical analysis which is to be done in the next section.
\section{Numerical Example and Analysis}
In this section we will present derive our control laws for a DC motor with the parameters given in \textbf{Table \ref{tab:DCMPARAM}}. 
\subsection{Speed Control}
\begin{figure}[H]
\centering
\includegraphics[scale=0.6]{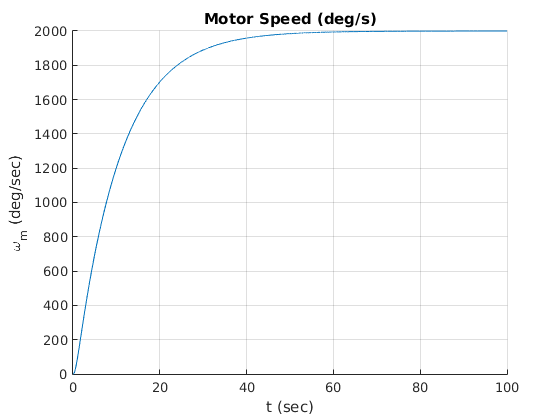}
\caption{Speed variation of the DC Motor parametrized in \textbf{Table \ref{tab:DCMPARAM}} under the control law defined by the output feedback gain in \eqref{PROJ_GAIN_SPEED}. The reference speed is $ \omega_r=2000\  \nicefrac{\mathrm{deg}}{\mathrm{sec}} $}. \label{SPEED_SIMULATION_NO_NOISE}
\end{figure}

\begin{figure}[H]
\centering
\includegraphics[scale=0.6]{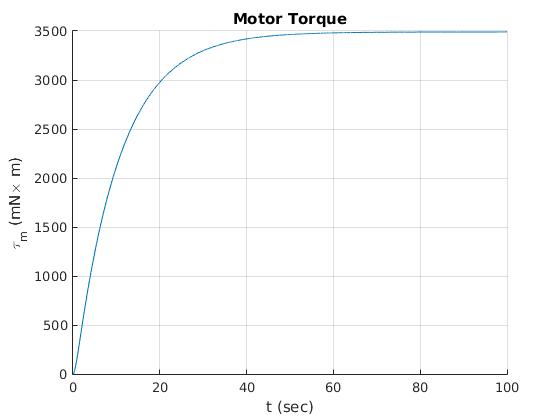}
\caption{Variation of the torque generated by the DC Motor parametrized in \textbf{Table \ref{tab:DCMPARAM}} under the control law defined by the output feedback gain in \eqref{PROJ_GAIN_SPEED}. The reference speed is $ \omega_r=2000\  \nicefrac{\mathrm{deg}}{\mathrm{sec}} $}. \label{SPEED_SIMULATION_NO_NOISE_TRQ}
\end{figure}

\begin{figure}[H]
\centering
\includegraphics[scale=0.6]{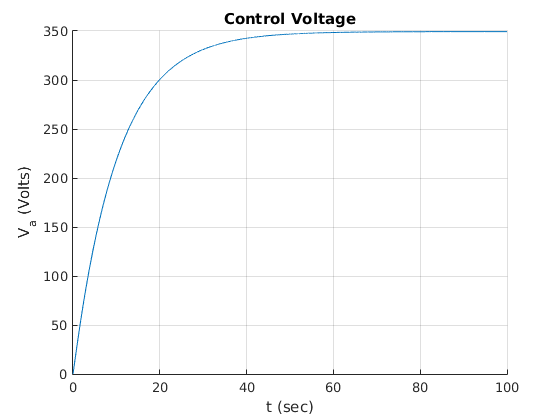}
\caption{Variation of the armature voltage required by the DC Motor parametrized in \textbf{Table \ref{tab:DCMPARAM}} under the control law defined by the output feedback gain in \eqref{PROJ_GAIN_SPEED}. The reference speed is $ \omega_r=2000\  \nicefrac{\mathrm{deg}}{\mathrm{sec}} $}. \label{SPEED_SIMULATION_NO_NOISE_VA}
\end{figure}

When the numerical values in \textbf{Table \ref{tab:DCMPARAM}} are substituted to \eqref{SPDCONTROL_SS}, one will obtain:
\begin{equation} \label{SPDCONTROL_SS_NUM}
\begin{bmatrix}
\dot{\epsilon}\\ \dot{\omega} \\ \dot{i_a}
\end{bmatrix}=\begin{bmatrix}
0 & 1 & 0 \\
0 & -10 & 1 \\
0 & -0.02 & -2 
\end{bmatrix}\begin{bmatrix}
\epsilon\\ \omega \\ i_a
\end{bmatrix}+\begin{bmatrix}
0 \\ 0 \\ 2
\end{bmatrix}V_a+\begin{bmatrix}
0 \\ 100 \\ 0
\end{bmatrix}\tau_L+\begin{bmatrix}
-1 \\ 0 \\ 0
\end{bmatrix}\omega_r
\end{equation}      
When one assumes that the disturbance torque equals to zero ($ \tau_L=0 $) and the reference velocity is a constant ($ \omega_r $ is a step function). The equation will be the same as \eqref{SPDCONTROL_SS_NO_TAU} with the system matrices same as in above. The full state feedback linear quadratic control can be obtained by invoking the matlab command \texttt{lqr(A,B,Q,R)} with \texttt{Q,R} being the matrices in the quadratic performance index shown in \eqref{QUADPI}. In this example they are taken as $ Q=qI_{3\times3} $ and $ R=1 $. When one invokes the MATLAB's \texttt{lqr} command for the given system in \eqref{SPDCONTROL_SS_NUM} with $ q=50 $, the resultant full state feedback control gain $ K_f $ in $ V_a=-K_fe $ (where $ e $ is defined in \textbf{Section \ref{SPD_CONTROL_SEC}}) is found as:
\begin{equation}
\label{Kf}
K_f=\begin{bmatrix}
7.071    &     0.903  &        6.204

\end{bmatrix}
\end{equation}
The above will yield the following closed loop spectrum:
\begin{subequations} \label{AcDV_SPEED}
\begin{equation} \label{Ac_SPEED} 
A_c=\begin{bmatrix}
         0  &   1.0000   &      0 \\
         0 & -10.0000  &  1.0000 \\
   -14.142   &   -1.8269   &   -14.409
\end{bmatrix}
\end{equation} 
\begin{equation} \label{D_SPEED} 
\Lambda_c=\begin{bmatrix}
    -0.098538     &       0     &      0\\
            0   &   -14.211      &      0\\
            0      &      0   &   -10.099

\end{bmatrix}
\end{equation}
\begin{equation} \label{V_SPEED} 
V_c=\begin{bmatrix}
     -0.71399  &  -0.016255  &  -0.098064 \\
     0.070355  &      0.231  &    0.99034 \\
      0.69662  &   -0.97282  &  -0.098014 
\end{bmatrix}
\end{equation}
\end{subequations}
where $ A_c=A-BK_f $, $ D_c=\lambda(A_c) $ and $ V_c $ is the eigenvectors corresponding to each element of $ D_c $ given in the order. From \eqref{SPDCONTROL_SS_NUM} one will recognize that the available state variables are the steady state errors of the integral of the velocity tracking error $ e_\epsilon(t)=\epsilon(t)-\epsilon(\infty) $ and the velocity itself $ e_\omega(t)=\omega(t)-\omega(\infty) $. This means that the output feedback matrix $ C $ is:
\begin{equation}
\label{C_MATR_MOTOR}
C=\begin{bmatrix}
1 & 0 & 0 \\
0 & 1 & 0
\end{bmatrix}
\end{equation} 
We should also note from the above equation that the number of available feedback lines is equal to 2 thus the number of eigenvalues that are to be retained from the closed loop full state feedback spectrum in \eqref{AcDV_SPEED} is also equal to 2. One has no control over the location of the third eigenvalue. In order to reduce the risk of an unstable mode, the desired retained eigenvalues among \eqref{D_SPEED} should be the two dominant ones in $ \Lambda_c $. Looking at \eqref{D_SPEED}, one can easily note that the dominant poles of the full state feedback closed loop are 
\begin{equation}
\label{DOM_CL_SPD}
\Lambda_r=\begin{bmatrix}
-0.098538   &      0  \\
0  & -10.099 
\end{bmatrix}
\end{equation}
and their corresponding eigenvectors are:
\begin{equation}
\label{DOM_CL_SPD_VR}
V_r=\begin{bmatrix}
-0.71399  &    -0.098064 \\
     0.070355  &    0.99034 \\
      0.69662  &  -0.098014 
\end{bmatrix}
\end{equation} 
So when one applies \eqref{PROJCONTROL_EQUATION}, the output feedback gain $ K_o $ will be found as:
\begin{equation}
\label{PROJ_GAIN_SPEED}
K_o=\begin{bmatrix}
0.89686  &   -0.32197
\end{bmatrix}
\end{equation} 
The above gain will yield an output feedback closed loop eigenvalues as shown below:
\begin{equation} \label{D_SPEED_OUFB} 
\Lambda_c=\begin{bmatrix}

    -0.098538      &      0   &         0\\
            0   &   -1.8025   &         0\\
            0   &         0   &   -10.099

\end{bmatrix}
\end{equation}
So only the $ 2^\mathrm{nd} $ eigenvalue is different from the full state feedback equivalent which is -14.211. So we obtained a stable output feedback based DC motor control system. In \textbf{Figures \ref{SPEED_SIMULATION_NO_NOISE},\ref{SPEED_SIMULATION_NO_NOISE_TRQ},\ref{SPEED_SIMULATION_NO_NOISE_VA}}, one can see the simulation results obtained when the control law with gain \eqref{PROJ_GAIN_SPEED} is applied. 

\begin{figure}[H]
\centering
\includegraphics[scale=0.6]{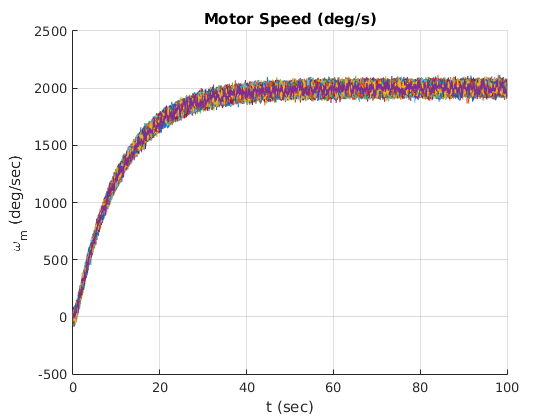}
\caption{Speed variation of the DC Motor parametrized in \textbf{Table \ref{tab:DCMPARAM}} under the control law defined by the output feedback gain in \eqref{PROJ_GAIN_SPEED}. The reference speed is $ \omega_r=2000\  \nicefrac{\mathrm{deg}}{\mathrm{sec}} $. In this simulation, disturbance torque is present as a Gaussian distributed random variable with mean $ \mu=0 $ and variance $ \sigma^2=200\ \mathrm{mN\cdot m} $}. \label{SPEED_SIMULATION_WTH_NOISE}
\end{figure}

\begin{figure}[H]
\centering
\includegraphics[scale=0.6]{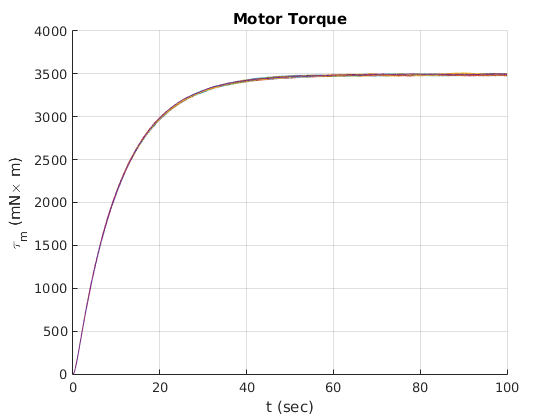}
\caption{Variation of the torque generated by the DC Motor parametrized in \textbf{Table \ref{tab:DCMPARAM}} under the control law defined by the output feedback gain in \eqref{PROJ_GAIN_SPEED}. The reference speed is $ \omega_r=2000\  \nicefrac{\mathrm{deg}}{\mathrm{sec}} $.In this simulation, disturbance torque is present as a Gaussian distributed random variable with mean $ \mu=0 $ and variance $ \sigma^2=200\ \mathrm{mN\cdot m} $}. \label{SPEED_SIMULATION_WTH_NOISE_TRQ}
\end{figure}

\begin{figure}[H]
\centering
\includegraphics[scale=0.6]{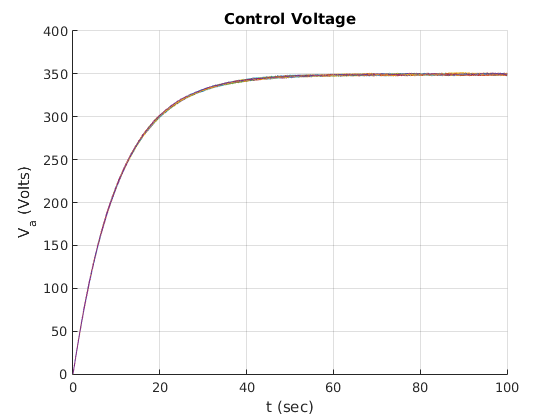}
\caption{Variation of the armature voltage required by the DC Motor parametrized in \textbf{Table \ref{tab:DCMPARAM}} under the control law defined by the output feedback gain in \eqref{PROJ_GAIN_SPEED}. The reference speed is $ \omega_r=2000\  \nicefrac{\mathrm{deg}}{\mathrm{sec}} $. In this simulation, disturbance torque is present as a Gaussian distributed random variable with mean $ \mu=0 $ and variance $ \sigma^2=200\ \mathrm{mN\cdot m} $}. \label{SPEED_SIMULATION_WTH_NOISE_VA}
\end{figure}

When one has a disturbance torque effective on the motor shaft $ (\tau_L) $, one will note the results presented in \textbf{Figures \ref{SPEED_SIMULATION_WTH_NOISE}, \ref{SPEED_SIMULATION_WTH_NOISE_TRQ}, \ref{SPEED_SIMULATION_WTH_NOISE_VA}}. Analysis of the results show that, the closed loop is working stable against disturbance torques. This might be seen as a violation of the \textbf{Theorem \ref{DSSTHEOREM}}. However, we have here to stress that \textbf{Theorem \ref{DSSTHEOREM}} is a sufficient not necessary condition. It is also a conservative inequality as it is transformed from Lyapunov equation \eqref{WDOTOPENFORM} by using upper/lower bound theorems \eqref{EQQUADBOUND}. Because of that, choosing \eqref{PROJ_GAIN_SPEED} which satisfies \textbf{Theorem \ref{DSSTHEOREM}} will guarantee disturbance-to-state stability but this does not mean that poles positioned near to $ j\omega $ axis will lead to instability under disturbance. One can see the simulation results in \textbf{Figures } when the eigenvalue nearest to $ j\omega $ axis is shifted to the position $ \lambda=-0.8 $. These results are obtained when $ K_o $ is replaced by:
\begin{equation}
\label{PROJ_GAIN_SPEED_CORR}
K_o=\begin{bmatrix}
4.4476  &   0.029499
\end{bmatrix}
\end{equation}         
The above will yield the closed loop eigenvalues as:
\begin{equation}
\label{Do_CORR}
\Lambda_o=\begin{bmatrix}
-0.8     &       0      &      0 \\
0    &   -1.101     &       0 \\
0     &       0   &   -10.099
\end{bmatrix}
\end{equation}
The second eigenvalue is not placed as expected but it does not violate \textbf{Theorem \ref{DSSTHEOREM}}. The evaluation of the gain in \eqref{PROJ_GAIN_SPEED_CORR} is performed by applying the orthogonal projection equation \eqref{PROJCONTROL_EQUATION} to the full state feedback spectrum obtained from a pole placement design which replaces the eigenvalue violating \textbf{Theorem \ref{DSSTHEOREM}} by a suitable one (i.e. $ \lambda=-0.8 $). However this will result in the loss of optimality provided by the linear quadratic regulator (suboptimal/near optimal controller). It should also be noted that, there is a very little change in the simulation results. 

The simulation based analysis of the disturbance torque effects are based on the repeated runs of the closed loop controlled model with a stochastic disturbance input. The disturbance torque $ \tau_L(t) $ is considered as a random variable with zero mean and a certain level of variance which is chosen to be less than the 10\% of the maximum value of the torque obtained. This exogenous input will effect the angular velocity of the motor and thus its position. In order to see the actual situation, the randomness of the disturbance will force one to repeat the simulations several times with the normally distributed disturbance torque active on the model. In this study, the number of repeats is chosen as $ N_{tst}=200 $ (200 times repeating). The numerical details of the disturbance torques for each group of simulation are either written in the figure captions or in the parts of the text referring to the illustrations.           

\begin{figure}[H]
\centering
\includegraphics[scale=0.6]{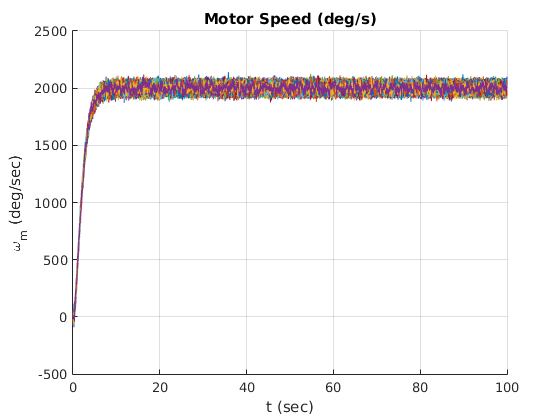}
\caption{Speed variation of the DC Motor parametrized in \textbf{Table \ref{tab:DCMPARAM}} under the control law defined by the output feedback gain in \eqref{PROJ_GAIN_SPEED}. The reference speed is $ \omega_r=2000\  \nicefrac{\mathrm{deg}}{\mathrm{sec}} $. In this simulation, disturbance torque is present as a Gaussian distributed random variable with mean $ \mu=0 $ and variance $ \sigma^2=200\ \mathrm{mN\cdot m} $. Here the control gain at \eqref{PROJ_GAIN_SPEED_CORR} is generating the control law.}. \label{SPEED_SIMULATION_WTH_NOISE_CORRP}
\end{figure} 
\begin{figure}[H]
\centering
\includegraphics[scale=0.6]{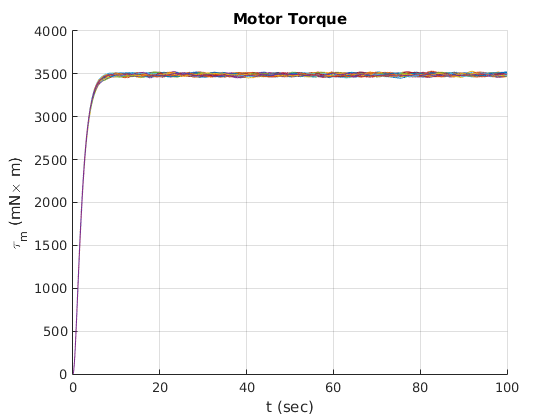}
\caption{Variation of the torque generated by the DC Motor parametrized in \textbf{Table \ref{tab:DCMPARAM}} under the control law defined by the output feedback gain in \eqref{PROJ_GAIN_SPEED}. The reference speed is $ \omega_r=2000\  \nicefrac{\mathrm{deg}}{\mathrm{sec}} $.In this simulation, disturbance torque is present as a Gaussian distributed random variable with mean $ \mu=0 $ and variance $ \sigma^2=200\ \mathrm{mN\cdot m} $. Here the control gain at \eqref{PROJ_GAIN_SPEED_CORR} is generating the control law.}. \label{SPEED_SIMULATION_WTH_NOISE_TRQ_CORRP}
\end{figure}
\begin{figure}[H]
\centering
\includegraphics[scale=0.6]{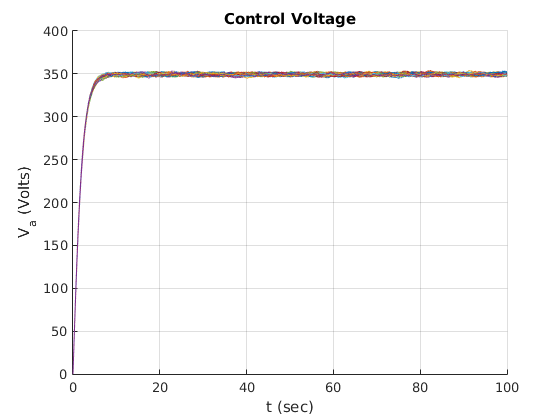}
\caption{Variation of the armature voltage required by the DC Motor parametrized in \textbf{Table \ref{tab:DCMPARAM}} under the control law defined by the output feedback gain in \eqref{PROJ_GAIN_SPEED}. The reference speed is $ \omega_r=2000\  \nicefrac{\mathrm{deg}}{\mathrm{sec}} $. In this simulation, disturbance torque is present as a Gaussian distributed random variable with mean $ \mu=0 $ and variance $ \sigma^2=200\ \mathrm{mN\cdot m} $. Here the control gain at \eqref{PROJ_GAIN_SPEED_CORR} is generating the control law.}. \label{SPEED_SIMULATION_WTH_NOISE_VA_CORRP}
\end{figure}

\subsection{Position Control}
In position control the numerics are mostly the same. Only some of the details on the state equations will differ. First of all, \eqref{SPDCONTROL_SS_NUM} will be replaced by:
\begin{equation} \label{PSNCONTROL_SS_NUM}
\begin{bmatrix}
\dot{e}_\theta \\ \dot{\omega} \\ \dot{i_a}
\end{bmatrix}=\begin{bmatrix}
0 & 1 & 0 \\
0 & -10 & 1 \\
0 & -0.02 & -2 
\end{bmatrix}\begin{bmatrix}
\epsilon\\ \omega \\ i_a
\end{bmatrix}+\begin{bmatrix}
0 \\ 0 \\ 2
\end{bmatrix}V_a+\begin{bmatrix}
0 \\ 100 \\ 0
\end{bmatrix}\tau_L
\end{equation}  
where $ e_\theta=\theta-\theta_r $ with $ \theta_r $ being the desired/reference position of the DC motor. The measured state variables in the above configuration are $ e_\theta $ and $ \omega $. Thus the output feedback mapping matrix is the same as that of \eqref{CMATR}. 
In addition, as the system matrices of \eqref{PSNCONTROL_SS_NUM} are numerically the same as that of \eqref{SPDCONTROL_SS_NUM}, the controller gains \eqref{Kf}, \eqref{PROJ_GAIN_SPEED}, \eqref{PROJ_GAIN_SPEED_CORR}, closed loop spectrum \eqref{AcDV_SPEED}, \eqref{DOM_CL_SPD}, \eqref{DOM_CL_SPD_VR}, \eqref{D_SPEED_OUFB} and \eqref{Do_CORR} will be same for position control problem provided that the quadratic performance coefficients are same $ Q=qI_{3\times 3} $ with $ q=50 $ and $ R=1 $. 

\begin{figure}[H]
\centering
\includegraphics[scale=0.6]{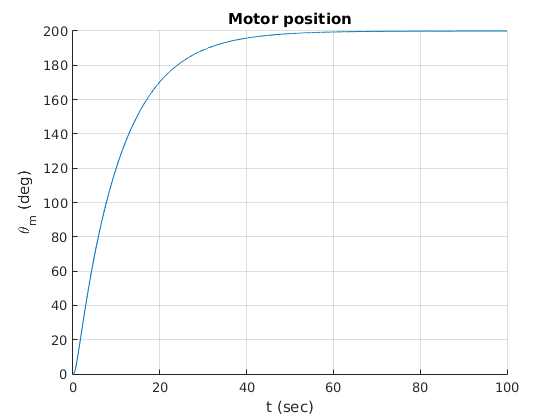}
\caption{Variation of the position of the DC Motor parametrized in \textbf{Table \ref{tab:DCMPARAM}} under the position control law defined by the output feedback gain in \eqref{PROJ_GAIN_SPEED} which is utilized as given in \eqref{PSN_CONTROL_LAWX}. The reference position is $ \theta_r=200\degree $}. \label{PSN_SIMULATION_NO_NOISE}
\end{figure}

\begin{figure}[H]
\centering
\includegraphics[scale=0.6]{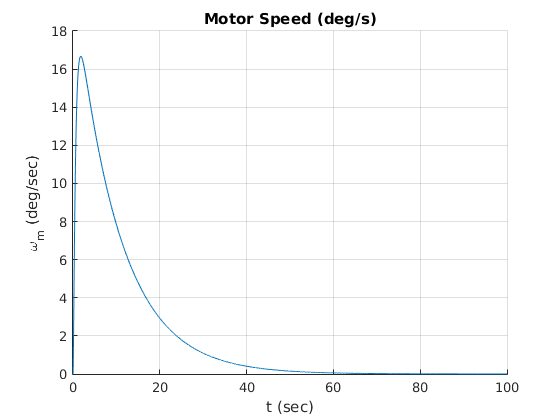}
\caption{Variation of the speed of the DC Motor parametrized in \textbf{Table \ref{tab:DCMPARAM}} under the position control law defined by the output feedback gain in \eqref{PROJ_GAIN_SPEED} which is utilized as given in \eqref{PSN_CONTROL_LAWX}. The reference position is $ \theta_r=200\degree $}. \label{PSN_SIMULATION_NO_NOISE_SPD}
\end{figure}

\begin{figure}[H]
\centering
\includegraphics[scale=0.6]{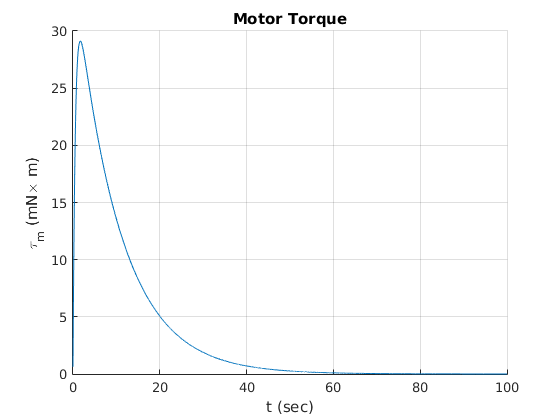}
\caption{Variation of the torque generated by the DC Motor parametrized in \textbf{Table \ref{tab:DCMPARAM}} under the position control law defined by the output feedback gain in \eqref{PROJ_GAIN_SPEED} which is utilized as given in \eqref{PSN_CONTROL_LAWX}. The reference position is $ \theta_r=200\degree $}. \label{PSN_SIMULATION_NO_NOISE_TRQ}
\end{figure}

\begin{figure}[H]
\centering
\includegraphics[scale=0.6]{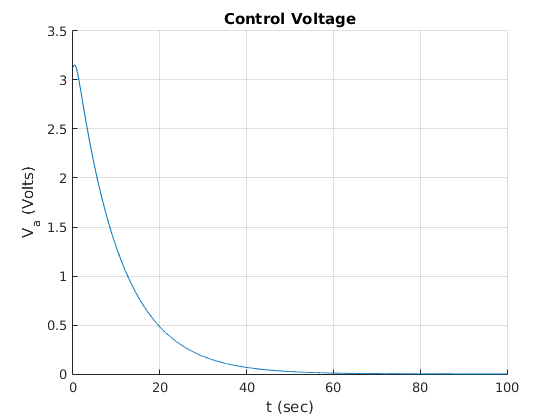}
\caption{Variation of the armature voltage required by the DC Motor parametrized in \textbf{Table \ref{tab:DCMPARAM}} under the position control law defined by the output feedback gain in \eqref{PROJ_GAIN_SPEED} which is utilized as given in \eqref{PSN_CONTROL_LAWX}. The reference position is $ \theta_r=200\degree $}. \label{PSN_SIMULATION_NO_NOISE_VA}
\end{figure}

\begin{figure}[H]
\centering
\includegraphics[scale=0.6]{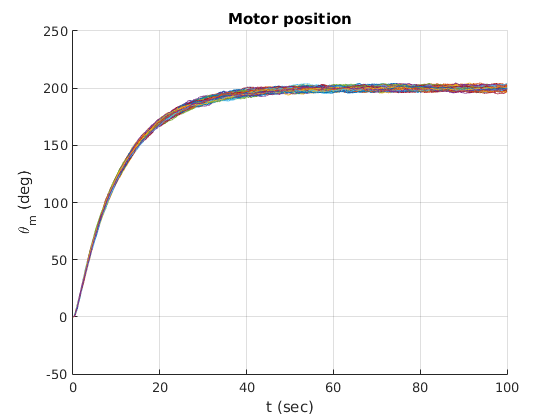}
\caption{Variation of the position of the DC Motor parametrized in \textbf{Table \ref{tab:DCMPARAM}} under the position control law defined by the output feedback gain in \eqref{PROJ_GAIN_SPEED} which is utilized as given in \eqref{PSN_CONTROL_LAWX}. The reference position is $ \theta_r=200\degree $. Here, the simulation is performed under the applied disturbance torques}. \label{PSN_SIMULATION_WTH_NOISE}
\end{figure}

\begin{figure}[H]
\centering
\includegraphics[scale=0.6]{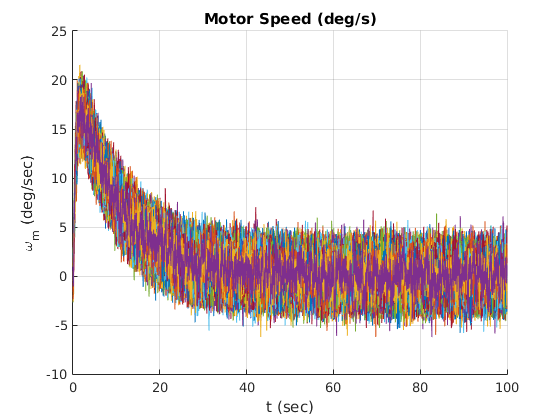}
\caption{Variation of the speed of the DC Motor parametrized in \textbf{Table \ref{tab:DCMPARAM}} under the position control law defined by the output feedback gain in \eqref{PROJ_GAIN_SPEED} which is utilized as given in \eqref{PSN_CONTROL_LAWX}. The reference position is $ \theta_r=200\degree $. Here, the simulation is performed under the applied disturbance torques}. \label{PSN_SIMULATION_WTH_NOISE_SPD}
\end{figure}

\begin{figure}[H]
\centering
\includegraphics[scale=0.6]{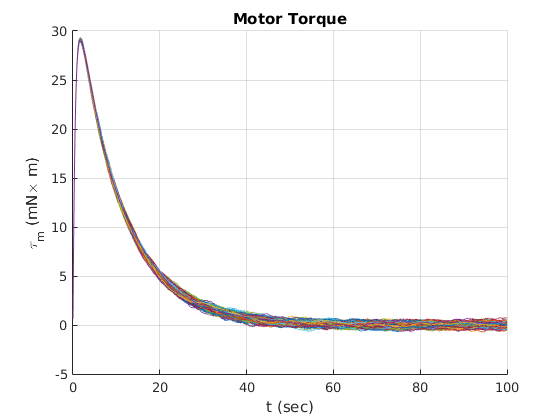}
\caption{Variation of the torque generated by the DC Motor parametrized in \textbf{Table \ref{tab:DCMPARAM}} under the position control law defined by the output feedback gain in \eqref{PROJ_GAIN_SPEED} which is utilized as given in \eqref{PSN_CONTROL_LAWX}. The reference position is $ \theta_r=200\degree $. Here, the simulation is performed under the applied disturbance torques}. \label{PSN_SIMULATION_WTH_NOISE_TRQ}
\end{figure}

\begin{figure}[H]
\centering
\includegraphics[scale=0.6]{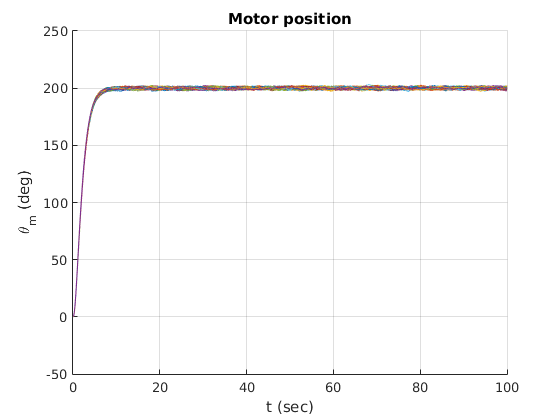}
\caption{Variation of the position of the DC Motor parametrized in \textbf{Table \ref{tab:DCMPARAM}} under the position control law defined by the output feedback gain in \eqref{PROJ_GAIN_SPEED} which is utilized as given in \eqref{PSN_CONTROL_LAWX}. The reference position is $ \theta_r=200\degree $. Here, the simulation is performed under the applied disturbance torques.  All the poles are satisfying \textbf{Theorem \ref{DSSTHEOREM}}.}. \label{PSN_SIMULATION_WTH_NOISE_CORR}
\end{figure}

\begin{figure}[H]
\centering
\includegraphics[scale=0.6]{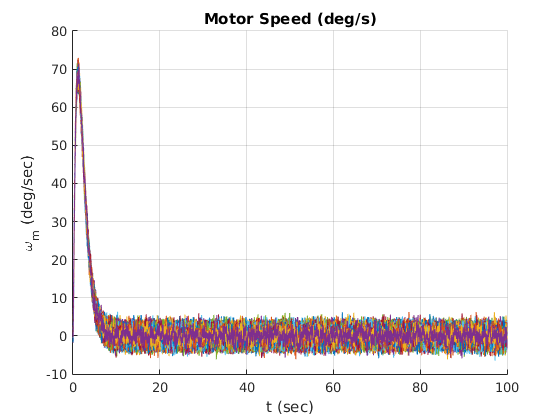}
\caption{Variation of the speed of the DC Motor parametrized in \textbf{Table \ref{tab:DCMPARAM}} under the position control law defined by the output feedback gain in \eqref{PROJ_GAIN_SPEED} which is utilized as given in \eqref{PSN_CONTROL_LAWX}. The reference position is $ \theta_r=200\degree $. Here, the simulation is performed under the applied disturbance torques. All the poles are satisfying \textbf{Theorem \ref{DSSTHEOREM}}.}. \label{PSN_SIMULATION_WTH_NOISE_SPD_CORR}
\end{figure}

\begin{figure}[H]
\centering
\includegraphics[scale=0.6]{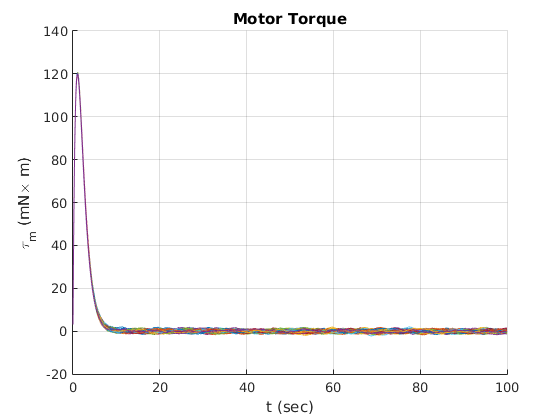}
\caption{Variation of the torque generated by the DC Motor parametrized in \textbf{Table \ref{tab:DCMPARAM}} under the position control law defined by the output feedback gain in \eqref{PROJ_GAIN_SPEED} which is utilized as given in \eqref{PSN_CONTROL_LAWX}. The reference position is $ \theta_r=200\degree $. Here, the simulation is performed under the applied disturbance torques. All the poles are satisfying \textbf{Theorem \ref{DSSTHEOREM}}.}. \label{PSN_SIMULATION_WTH_NOISE_TRQ_CORR}
\end{figure}

\begin{figure}[H]
\centering
\includegraphics[scale=0.6]{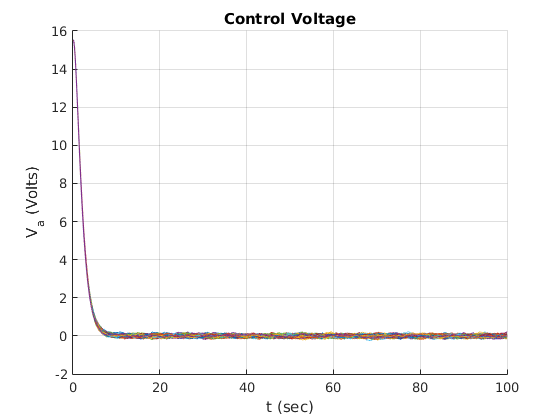}
\caption{Variation of the armature voltage required by the DC Motor parametrized in \textbf{Table \ref{tab:DCMPARAM}} under the position control law defined by the output feedback gain in \eqref{PROJ_GAIN_SPEED} which is utilized as given in \eqref{PSN_CONTROL_LAWX}. The reference position is $ \theta_r=200\degree $. Here, the simulation is performed under the applied disturbance torques. All the poles are satisfying \textbf{Theorem \ref{DSSTHEOREM}}.}. \label{PSN_SIMULATION_WTH_NOISE_VA_CORR}
\end{figure}

In \textbf{Figures \ref{PSN_SIMULATION_NO_NOISE},\ref{PSN_SIMULATION_NO_NOISE_TRQ},\ref{PSN_SIMULATION_NO_NOISE_VA}}, one will be able to see the position tracking simulation under noise free environment when the control law defined by the gain $ K_o $ in \eqref{PROJ_GAIN_SPEED} is applied as:
\begin{equation}
V_a=-K_o\begin{bmatrix}e_\theta & \omega\end{bmatrix}^T
\label{PSN_CONTROL_LAWX}
\end{equation}
Using the same configuration that resulted \textbf{Figures \ref{PSN_SIMULATION_NO_NOISE},\ref{PSN_SIMULATION_NO_NOISE_SPD},\ref{PSN_SIMULATION_NO_NOISE_TRQ},\ref{PSN_SIMULATION_NO_NOISE_VA}} the simulation in a noisy environment (with disturbance torque) yields the results shown in \textbf{Figures \ref{PSN_SIMULATION_WTH_NOISE},\ref{PSN_SIMULATION_WTH_NOISE_SPD},\ref{PSN_SIMULATION_WTH_NOISE_TRQ},\ref{PSN_SIMULATION_NO_NOISE_VA}}. In this case, a disturbance torque is effective on the motor shaft and it is modeled by a Gaussian distributed source with zero mean and variance $ \sigma^2=0.01\  \mathrm{N}\cdot\mathrm{m}$. As we have done in the case of speed control, we will present the results of the simulation when the smallest eigenvalue at $ \lambda=-0.098538 $ is moved to $ \lambda=-0.8 $ in \textbf{Figures \ref{PSN_SIMULATION_WTH_NOISE_CORR},\ref{PSN_SIMULATION_WTH_NOISE_SPD_CORR},\ref{PSN_SIMULATION_WTH_NOISE_TRQ_CORR},\ref{PSN_SIMULATION_WTH_NOISE_VA_CORR}}.

The examples given in this section are to demonstrate the approaches presented in \textbf{Section \ref{LQPC}} which is the linear quadratic projective control approach. The purpose is to present the methodology such that, interested readers can replicate the procedure. Thus, given a single reference position or speed (final target position/speed) we presented the simulation results. For testing our controllers under noise due to the disturbance torques we presented repeated trials where the disturbance torque appears as a normally distributed random variable. In each run the disturbance profile will be different due to its randomness so one can reflect those analyses as Monte Carlo methods \cite{caflisch1998monte} which rely on repeated samples of random data to obtain the performance of an algorithm when there are parameters with uncertainty. With this view, one can also treat this approach as a robust stability test. Nevertheless, the theoretical stability discussion (\textbf{Theorem \ref{DSSTHEOREM}}) a better approach which is considered a general methodology regardless of the type of the disturbance torques.       

\begin{figure}[H]
\centering
\includegraphics[scale=0.6]{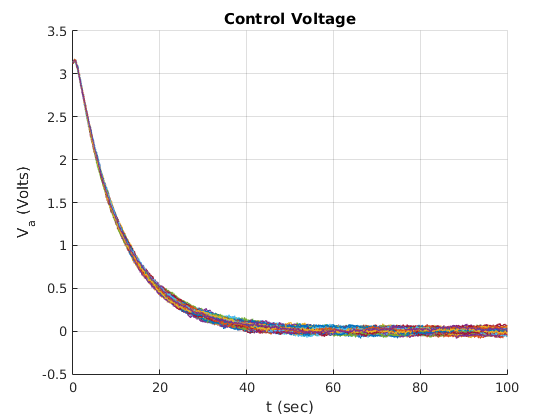}
\caption{Variation of the armature voltage required by the DC Motor parametrized in \textbf{Table \ref{tab:DCMPARAM}} under the position control law defined by the output feedback gain in \eqref{PROJ_GAIN_SPEED} which is utilized as given in \eqref{PSN_CONTROL_LAWX}. The reference position is $ \theta_r=200\degree $. Here, the simulation is performed under the applied disturbance torques}. \label{PSN_SIMULATION_WTH_NOISE_VA}
\end{figure}

% \section{Authors' contact details, photos and biographies}
%Authors' biographies, photos and contact details will be submitted separately to the editorial office in case the
%article is accepted for publication, together with the final paper version. They will be inserted at the end of the
%article by the editorial office (for information how this looks like see published articles from either vol.~50, no.~3-4
%or any later number).

\section{Conclusion}\label{sec:concl}
In this work, we presented a direct current electrical motor control by linear quadratic projective control. The chosen methodology helps the designer to eliminate the feedback from armature current which will increase cost of feedback and incorporate high noise due to its amplification in the signal conditioning circuitry. The simulations reveal that under both ideal conditions and noisy environment (due to the disturbance torques on the motor shaft), the controller can handle its operation very well and works stable. In addition to simulations, a theoretical discussion on the stability of the closed loop against the disturbance torques is given which is based on the input-to-state stability concept. The theoretical development is fairly conservative as it is developed from the conversion of the equations related to the Lyapunov's second method to inequalities through its manipulation by upper and lower bound lemmas. This result is also seen from the simulations. The designs which does not satisfy the disturbance to state stability theorem (\textbf{Theorem \ref{DSSTHEOREM}}) can handle the noises without going into instability. However, one can not simulate all kinds of disturbance torques as disturbances are a particular group of stochastic processes as their name implies. So a design satisfying the disturbance to state stability theorem is expected to guarantee the closed loop's stability against large disturbance torques. In addition to that, in both speed and position controllers the larger eigenvalues provide a closed loop with faster response times. 

A future study based on this work can be the application of different linear control techniques to the same problem and repeat the theoretical and numerical analysis performed in this work on them.

% use /balance at an appropriate place before the paper end to get approximately equally long columns on the last text
% page
\balance

% Start appendix with \appendixAutomatika command
\appendixAutomatika
% \section{Appendix section}\label{sec:appendixa}

%This is an example of appendix section. This is an example of appendix section. This is an example of appendix section.
%This is an example of appendix section.

% \subsection{Appendix subsection}\label{subsec:appendix_subsec}

%This is an example of appendix subsection. This is an example of appendix subsection. This is an example of appendix
%subsection. This is an example of appendix subsection. This is an example of appendix subsection. This is an example of
%appendix subsection. This is an example of appendix subsection. This is an example of appendix subsection. This is an
%example of appendix subsection. This is an example of appendix subsection. This is an example of appendix subsection. 
%This is an example of appendix subsection. This is an example of appendix subsection. 
\appendixendAutomatika
% End appendix with \appendixendAutomatika command

%% Use \section* for Acknowledgment
%\section*{Acknowledgment}
%Type your acknowledgment and thanks here. Note that it is customary to put acknowledgment only in the accepted, final
%paper version and not to put it in the version submitted for review.

%=============
% BIBLIOGRAPHY
%=============
% The authors are encouraged to use a bibliography generated by BibTeX as a .bbl file, instead writing it manually. The
% bibliography style to be used is the "ieeetr" style. See the example.bib for reference
\small
\bibliographystyle{ieeetr}
\bibliography{AutomatikaTemplateLatex.bib}
\balance
\end{document}